\documentclass[11pt,a4]{article}
\newtheorem{theorem}{Theorem}%[section]
\newtheorem{proposition}[theorem]{Proposition}

     \newtheorem{propositionf}[theorem]{``–½'è"}

\begin{document}

\begin{center}
{\Large
\bf Elementary divisors of Cartan matrices for symmetric groups}

\vspace{5mm}

Katsuhiro UNO and Hiro-Fumi YAMADA
\end{center}

\vspace{3mm}

\section{Introduction}

       The purpose of this note is to give a simple expression of the
elementary divisors of the Cartan matrices for the symmetric
groups.

       Let $G$ be a finite group and $k$ be an algebraically closed field
of characteristic $p > 0$.  The group algebra $kG$, which affords the
left regular representation of $G$, is a direct sum of indecomposable representations.
There is a natural one-to-one correspondence between the equivalence
classes of indecomposable summands of $kG$ and those of irreducible
representations.  Let $c_{\lambda \mu}$ be the multiplicity of
the irreducible representation $F_{\mu}$ occurring as a composition factor
of the indecomposable summand $U_{\lambda}$ of $kG$.
The Cartan matrix is, by definition $C = (c_{\lambda \mu})$, and is an
$\ell \times \ell$ square matrix, where $\ell$ is the number of inequivalent
irreducible representations.
The following is known for the elementary divisors of $C$.
Let  $\{x_1, \cdots, x_{\ell}\}$ be a set of representatives of $p$-regular conjugacy classes of  $G$.  Then the elementary divisors
of $C$ are given by $\{|Z_G(x_i)|_p ; i = 1, \cdots, \ell\}$,
where $Z_G(x)$ denotes the centralizer of $x$ in $G$ \cite{NT}.
In the case of the symmetric groups, explicit computations are made in \cite{O1}  (see also \cite{BO1}), by using generating functions.

      Here we give in this note a simple expression of the elementary divisors
for the symmetric groups by utilizing the one-to-one correspondence between
the $p$-regular partitions and $p$-class regular partitions. By our result, one can obtain them by looking at $p$-regular partitions only.
      The case $p=2$ is of special interest.  We can state the block version
of our expression.  As is well-known, the elementary divisors of the Cartan
matrix depend only on the weight $w$ of the block.  For the case $p=2$ we
determine the elementary divisors of the given weight by looking at a certain
abacus.  It is interesting that the same abacus plays a role in the
$Q$-function realization of the basic representation of the affine Lie
algebra $A^{(1)}_1$ \cite{NY}.
%\vspace{5mm}

\section{The Glaisher map and elementary divisors}
%\vspace{5mm}

      In this section $p$ always denotes a fixed prime number.
A partition $\lambda$ is said to be $p$-regular if it does not contain
$p$ equal parts.  We denote the set of all $p$-regular partitions of $n$
by $P^{r(p)}(n)$.  A partition $\lambda$ is said to be $p$-class regular
if $p$ does not divide any part of $\lambda$.  We denote the set of
all $p$-class regular partitions of $n$ by $P^{c(p)}(n)$.  There is a
bijection between $P^{r(p)}(n)$ and $P^{c(p)}(n)$, which is defined as follows.
Let $\lambda=(\lambda_1, \cdots, \lambda_{\ell})$ be a $p$-regular
partition.  Write each part as $\lambda_i = p^{a_i}q_i$ with
$(p,q_i) = 1$.  Let $\mu(i)$ be the rectangular partition of $\lambda_i$
given by $\mu(i) = (q_i, \cdots, q_i)$ with length $p^{a_i}$.
Suppose that $q_{j_1}\geq \cdots \geq q_{j_{\ell}}$.  Then let
$\tilde{\lambda}$ be the vertical concatenation
$(\mu(j_1), \cdots, \mu(j_{\ell}))$, which is $p$-class regular.
For example, if $p=2$ and $\lambda =(5,4,3)$, then
$\tilde{\lambda} =(5,3,1,1,1,1)$. Then, the map $\gamma :P^{r(p)}(n)  \to P^{c(p)}(n)$ sending any $\lambda$ in $P^{r(p)}(n)$ into $\tilde{\lambda}$ is well defined, and it is
 easily seen that $\gamma$ gives a bijection. This $\gamma$ is sometimes called the Glaisher map. Note that $\ell(\tilde{\lambda}) \geq
\ell(\lambda)$. The elementary divisors of the $p$-Cartan matrix can be described as follows.

%\vspace{5mm}

\begin{theorem}\label{ele} The elementary divisors of the $p$-Cartan matrix for
the symmetric group $S_n$ on $n$ letters are given by
$$\{p^{\frac{\ell(\tilde{\lambda})-\ell(\lambda)}{p-1}} ; \lambda \in
P^{r(p)}(n) \}.$$
\end{theorem}
\vspace{5mm}

\noindent Proof.  Take $\tilde{\lambda} \in P^{c(p)}(n)$ and fix a representative
$x_{\tilde{\lambda}}$ of the corresponding $p$-regular class in $S_n$.
We will prove that
$$|Z_{S_n}(x_{\tilde{\lambda}})|_p =
p^{\frac{\ell(\tilde{\lambda})-\ell(\lambda)}{p-1}}.$$
Putting $m_k$ to be the multiplicity of the natural number $k$
as the parts of $\tilde{\lambda}$, we  write
$\tilde{\lambda}$ exponentially as $(1^{m_1}, \cdots, k^{m_k}, \cdots)$.
Let $m_k = \sum_{i \geq 0} b_i^{(k)}p^i$ be the
$p$-adic expansion of $m_k$.  ($0\leq b_i^{(k)} <p $.) Since
$$|Z_{S_n}(x_{\tilde{\lambda}})| = \prod_k k^{m_k}(m_k!),$$
we have
$$|Z_{S_n}(x_{\tilde{\lambda}})|_p = \prod_k (m_k!)_p = p^{\Sigma_k
d_p(m_k)},$$
where
$$d_p(m) = \sum_{i \geq 1}[\frac{m}{p^i}].$$
A simple computation shows that
$$d_p(m_k) = \frac{1}{p-1}(m_k - \sum_{i \geq 0} b_i^{(k)})$$
and we have
$$\sum_k d_p(m_k) = \frac{1}{p-1}(\sum_k m_k - \sum_{k,i} b_i^{(k)})
= \frac{1}{p-1}(\ell(\tilde{\lambda})-\ell(\lambda)).$$
Thus, the assertion holds.

\vspace{5mm}
The invariant $(\ell (\tilde{\lambda }) - \ell (\lambda ))/(p-1)$ can also be interpreted as follows. As in the paragraph preceding Theorem 1, for $\lambda = (\lambda _1, \lambda _2, \cdots , \lambda _{\ell})$ we write   $\lambda_i = p^{a_i}q_i$ with
$(p,q_i) = 1$.  We then have
$$\frac{1}{p-1}(\ell(\tilde{\lambda})-\ell(\lambda)) = \sum _i \frac{p^{a_i}-1}{p-1}  = \sum _{i, \ a_i \ne 0} (1 + p + \cdots + p^{a_i-1}).$$
Let $p\lambda$ denote the partition $(p\lambda _1, p\lambda _2, \cdots , p\lambda _{\ell})$. Then, it is easy to see that
$$\frac{1}{p-1}(\ell(\tilde{p\lambda})-\ell(p\lambda)) = \sum _i \frac{p^{a_i+1}-1}{p-1}   = \sum _{i} (1 + p + \cdots + p^{a_i}) = e(\lambda ).$$
Here, the invariant $e(\lambda)$ is the one first introduced incorrectly on p.54 in \cite{O1} and then correctly on p.299 in \cite{BO1}. In these articles, the multiplicity of $p^i$ as an elementary divisor of the Cartan matrix is given in terms of the number of $p$-regular partition $\mu $ of $k$ with $e(\mu ) = i$ and the invariant $m(w-k)$. (Proposition 3.19 of \cite{O1}) Here $w$ is the weight of the block. However, we know only the generating function of $m$ and to obtain the multiplicity it is necessary to know the numbers of partitions $\mu$ of all $k$ smaller than or equal to $w$.  Thus it is quite interesting that from Theorem 1, in order to compute   the elementary divisor corresponding to $\tilde{\lambda}$, it suffices to look at the partitions $\lambda$ and $\tilde{\lambda}$ only. On the other hand, it should be noticed  that the arguments in \cite{O1} and \cite{BO1} are given block wise, while ours is stated for the entire group. The block  wise statement in our setting can be done only for $p=2$ at present and is given in the next section.

\section{A block version for $p=2$}

%\vspace{5mm}

     In this section we restrict our attention to the special case $p=2$.
The Glaisher map turns out to the bijection between the set of all strict
partitions and  that of all odd partitions.
Theorem \ref{ele} tells us that elementary divisors have the form  $2^{\ell(\tilde{\lambda})-\ell(\lambda)}$ for  strict partitions $\lambda$.

We introduce the following "H-abacus"(cf. [NY]).
They represent strict partitions. For example,
The H-abacus of $\lambda = (9,5,3,2)$ is shown below.

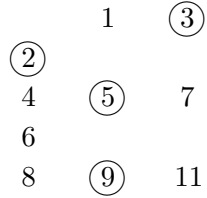
\begin{center}
\begin{picture}(80,80)
\put (40,65){1}
\put (70,65){3}
\put (10,50){2}
\put (10,35){4}
\put (40,35){5}
\put (70,35){7}
\put (10,20){6}
\put (10,5){8}
\put (40,5){9}
\put (68,5){11}
\put (12.5,53){\circle{13}}
\put (72.5,68){\circle{13}}
\put (42.5,38){\circle{13}}
\put (42.5,8){\circle{13}}
\end{picture}
\end{center}

For a strict partition we put a set of beads on the assigned
positions.  Any two beads do not pile up.  From the H-abacus of
the given strict partition $\lambda$, we obtain the
"H-core" $\lambda^H$ by moving and removing the beads as follows,
\vspace{5mm}

(1)  Move a bead one position up along the leftmost runner.

(2)  Remove a bead at the position 2.

(3)  Move a bead one position up along the runner of 1 or of 3.

(4)  Remove the two beads at the positions 1 and 3 simultaneously.
\vspace{5mm}

\noindent
The H-cores are thus characterized by the "stalemates", which
constitute the set
$$HC = \{\phi, (4m+1, 4m-3, \cdots, 5,1), (4m+3, 4m-1, \cdots, 7,3)
\quad (m \geq 0)\}.$$
For example, the H-core of the above $\lambda=(9,5,3,2)$ is
$\lambda^H = (1)$.

\vskip 0.3cm

Here, notice that the number of nodes in every H-core is a triangular number, $m(m+1)/2$, and conversely, for any triangular number $r$, there is a unique H-core with $r$ nodes. Recall that a similar thing is true for 2-cores since each 2-core has the form $\Delta_m = (m,m-1, \cdots, 2,1)$. Hence, there is a unique bijection between $HC$ and the set of 2-cores that preserves the number of nodes. In fact, the bijection can be obtained by applying "unfolding", which is defined as taking the hook lengths
of the main diagonal in the Young diagram. Namely, we have
$$HC = \{\Delta_m^{u}; m \geq 0\}, $$
where  $\lambda^{u}$
stands for the "unfolding" of $\lambda$.  For example,
 $\Delta_4^{u} = (7,3)$. For detail, see \S 4 of \cite{NY}.
\vskip 0.3cm

In the computation of the elementary divisors block wise, the H-core $\lambda^H$
of $\lambda$ plays an essential role, instead of the 2-core $\lambda^c$.
Let $\lambda$ be a strict partition of $n$.  To obtain the H-core
$\lambda^H$, we perform the moving (or removing) of the beads on the
H-abacus.  Let each moving (1) and (2) have H-weight 1, while (3) and (4)
have H-weight 2.  Thus $\lambda$ has its own H-weight,
which is, of course, equal to $w=(|\lambda|-|\lambda^H|)/2$.
For example $\lambda=(9,5,3,2)$ has H-weight $w=9$.

Let $SP(n)_w$ (resp. $SP(n)^w$) be the set of those strict partitions of $n$
with 2-weight $w$ (resp. H-weight $w$).  For example,
$SP(7)_2 =\{(6,1), (4,3)\}$, while $SP(7)^2 = \{(7), (4,3)\}$.

\vskip 0.3cm

The number of strict partitions of $n$ with a given H-core can be computed as follows.   Let $p(n)$ be the number of partitions of $n$, and let $q(n)$ be the number of strict partitions of $n$. Recall that $q(n)$ is also the number of partitions of $n$ with only odd parts.
For a partition $\lambda$, let $\lambda ^{(2')}$ and $\lambda ^{(2)}$ denote the partition obtained from $\lambda$ by picking up all the odd and even parts, respectively. Then, we have $\lambda ^{(2)} = 2 \mu = (2\mu _1, 2\mu _2, \cdots , 2 \mu_{\ell})$ for some partition $\mu$. The partitions $\lambda ^{(2')}$ and $\mu$ are determined uniquely by $\lambda$, and conversely, $\lambda$ is determined by an odd partition $\nu$ and a partition $\mu$ with $\lambda ^{(2')} = \nu$ and $\lambda ^{(2)} = 2 \mu$. This gives us
$$p(n) = \sum_{k =0}^{n} q(k)p(\frac{n-k}{2}),$$
where $p(x) =0$  if $x$ is not an integer.

\vskip 0.3cm
Note that the H-core of $\lambda$ depends only on $\lambda ^{(2')}$. Moreover, by arguments similar to those  in I.4 of \cite{O2}, partitions $\lambda$ in $SP(n)^w$ is determined uniquely by a pair of a strict partition $\mu$ with $2\mu = \lambda ^{(2)}$ and a partition $\nu$ obtained as the Frobenius symbol given by the 4-bar quotients, denoted by $\lambda _1$ and $\lambda _3$ in the notation in \cite{O2}, such that $|\mu | +2|\nu | = w$. Though some arguments in I.4 in \cite{O2} are given for odd numbers $p$ rather than 4, the above conclusion can be proved equally well for $p=4$. This shows that $|SP(n)^w| = \sum_{k =0}^{w} q(k)p(\frac{w-k}{2})$. Hence, we have $|SP(n)^w| = p(w)$.
\vskip 0.3cm

Recall also that for a given triangular number $r$, there are unique 2-core $\Delta$ and H-core $\Delta ^{u}$ with $r$ nodes.  Thus the set $SP(n)_w$ and $SP(n)^w$ consist of  strict partitions $\lambda$ of $n$ with $\lambda ^c = \Delta$ and  $\lambda ^H = \Delta ^{u}$, respectively, with $n-2w = r$.  It is known that the number of strict partitions of $n$ with a given 2-core depends only on the weight $w$ (Theorem 1 of \cite{O}) and in fact we have $|SP(n)_w| = p(w)$. Thus,  the following holds.

%\vspace{5mm}
\begin{proposition}  We have $|SP(n)_w| = |SP(n)^w| = p(w)$.
\end{proposition}
%\vspace{5mm}

By refining the above arguments, the  multiplicity of  $2^i$ as an elementary divisor of the Cartan matrix can be computed.

First note that $\ell(\tilde{\lambda})-\ell(\lambda) = \ell(\widetilde{\lambda ^{(2)}})-\ell(\lambda ^{(2)})$. In particular, $\ell (\tilde{\lambda} )-\ell (\lambda)$ depends only  on $\lambda ^{(2)}$. In other words, using the notation preceding the proposition,  $\ell(\tilde{\lambda})-\ell(\lambda)$ is determined only by $\mu$ and not by $\nu$. Hence the
  number of strict partitions $\lambda$ of $n$ with H-weight $w$ and
$\ell(\tilde{\lambda})-\ell(\lambda) = i$ is
$$\sum_{k =0}^{w} p_0^i(k)p(\frac{w-k}{2}),$$
where $p_0^i(k)$ is the number of $\mu \in P^{r(2)}(k)$ with $e(\mu ) =\ell (\tilde{2\mu} )-\ell (2\mu) = i$. Since the function $m(n)$ defined on p.53 in \cite{O1} satisfies
$$m(n) = p(\frac{n}{2}),$$
(see (3.12) in \cite{O1}), the formula  $\sum_{k =0}^{w} p_0^i(k)m(w-k)$ in Proposition 3.19 in \cite{O1} implies the following.

\begin{theorem}  Let $B$ be the unique 2-block of the symmetric group on $n$ letters corresponding to the 2-core $\Delta _m$. Then, the elementary divisors of the Cartan matrix for $B$ is given by
$$\{ \ 2 ^{\ell (\tilde{\lambda} )-\ell (\lambda)} \ | \ \lambda \in P^{r(2)}(n), \  \lambda ^H = \Delta _m^{u} \ \}.$$

In particular, the multiplicity of $2^i$ as an elementary divisor of the Cartan matrix depends only on the H-weight, which is equal to the 2-weight of $B$.
\end{theorem}

The above theorem shows that H-cores and H-weights give a combinatorial meaning for the formula in Proposition 3.19 in \cite{O1} for $p=2$. For odd primes, we do not know such an explanation.
\vspace{5mm}

\noindent Example.  For $n=7$, there are five strict partitions described below.
$$\begin{array}{llcll}
\lambda & \tilde{\lambda } & 2^{\ell (\tilde{\lambda} )-\ell (\lambda)} & \lambda ^c & \lambda ^H \\
\hline
(7) & (7) & 2^0& \Delta _1 = (1) &  \Delta _3^u = (3) \\
(6,1) & (3,3,1) & 2^1 & \Delta _3 =(2,1) & \Delta _1 ^u=(1)\\
(5,2 ) & (5,1,1) & 2^1 &\Delta _1 =(1) & \Delta _1^u= (1)\\
(4,3) & (3,1,1,1,1) & 2^3 &\Delta _3 =(2,1) &  \Delta _3 ^u = (3)\\
(4,2,1) & (1,1,1,1,1,1,1) & 2^4 &\Delta _1 =(1) &  \Delta _1^u=(1)\\
\end{array}$$

Thus, the elementary divisors of the Cartan matrix of the principal 2-block corresponding to $ \Delta _1$ is $\{ \ 2^1, 2^1, 2^4 \ \}$, and for the 2-block corresponding to  $\Delta _3$ we have $\{ \ 2^0, 2^3 \ \}$. However, irreducible 2-modular representations of the principal 2-block are labeled by $\{  (7), (5,2), (4,2,1)  \}$, and those of the non-principal 2-block by $\{  (6,1), (4,3)  \}$. See \cite{JK}.

\end{document}